\documentstyle[11pt]{article}
\title{An alternative proof of the extended Saalsch\"utz summation theorem 
for the ${}_{r+3}F_{r+2}(1)$ series with applications}
\author{
Y. S. Kim,\footnote{Department of Mathematics Education, Wonkwang University, Iksan, Korea
E-mail: yspkim@wonkwang.ac.kr}
\ \ Arjun. K. Rathie\footnote{Department of Mathematics, Central University of Kerala,
Kasaragad 671328, Kerala, India.
E-Mail: akrathie@cukerala.edu.in} \ \ and 
\ R. B. Paris\footnote{School of Computing, Engineering and Applied Mathematics, University of Abertay Dundee, Dundee DD1 1HG, UK.
E-Mail: r.paris@abertay.ac.uk}\ \footnote{Corresponding author}
 \\}
\begin{document}
\def\f#1#2{\mbox{${\textstyle \frac{#1}{#2}}$}}
\def\dfrac#1#2{\displaystyle{\frac{#1}{#2}}}
\def\boldal{\mbox{\boldmath $\alpha$}}
\newcommand{\bee}{\begin{equation}}
\newcommand{\ee}{\end{equation}}
\newcommand{\lam}{\lambda}
\newcommand{\ka}{\kappa}
\newcommand{\al}{\alpha}
\newcommand{\om}{\omega}
\newcommand{\Om}{\Omega}
\newcommand{\fr}{\frac{1}{2}}
\newcommand{\fs}{\f{1}{2}}
\newcommand{\g}{\Gamma}
\newcommand{\br}{\biggr}
\newcommand{\bl}{\biggl}
\newcommand{\ra}{\rightarrow}
\newcommand{\mbint}{\frac{1}{2\pi i}\int_{c-\infty i}^{c+\infty i}}
\newcommand{\mbcint}{\frac{1}{2\pi i}\int_C}
\newcommand{\mboint}{\frac{1}{2\pi i}\int_{-\infty i}^{\infty i}}
\newcommand{\gtwid}{\raisebox{-.8ex}{\mbox{$\stackrel{\textstyle >}{\sim}$}}}
\newcommand{\ltwid}{\raisebox{-.8ex}{\mbox{$\stackrel{\textstyle <}{\sim}$}}}
\renewcommand{\topfraction}{0.9}
\renewcommand{\bottomfraction}{0.9}
\renewcommand{\textfraction}{0.05}
\newcommand{\mcol}{\multicolumn}
\date{}
\maketitle
\begin{abstract}
A simple proof of a new summation formula for a terminating ${}_{r+3}F_{r+2}(1)$ hypergeometric series, representing an extension of Saalsch\"utz's formula for a ${}_3F_2(1)$ series, is given for the case of $r$ pairs of numeratorial and denominatorial parameters differing by positive integers. Two applications of this extended summation theorem are discussed. The first application extends two identities given by Ramanujan
and the second, which also employs a similar extension of the Vandermonde-Chu summation theorem for the ${}_2F_1$ series, extends certain reduction formulas for the Kamp\' e de F\'eriet function of two variables given by Exton and Cvijovi\'c and Miller. 

\vspace{0.4cm}

\vspace{0.3cm}

\end{abstract}

\vspace{0.3cm}

\begin{center}
{\bf 1. \  Introduction}
\end{center}
\setcounter{section}{1}
\setcounter{equation}{0}
\renewcommand{\theequation}{\arabic{section}.\arabic{equation}}
The generalized hypergeometric function ${}_pF_q(x)$ is defined for complex parameters and argument by the series \cite[p.~40]{S}
\bee\label{e11}
{}_pF_q\left[\!\!\begin{array}{c} a_1, a_2, \ldots ,a_p\\b_1, b_2, \ldots ,b_q\end{array}\!; z\right]=
\sum_{k=0}^\infty \frac{(a_1)_k (a_2)_k \ldots (a_p)_k}{(b_1)_k (b_2)_k \ldots (b_q)_k}\,\frac{z^k}{k!}.
\ee
When $q\geq p$ this series converges for $|z|<\infty$, but when $q=p-1$ convergence occurs when $|z|<1$.
However, when only one of the numeratorial parameters $a_j$ is a negative integer or zero,
then the series always converges since it is simply a polynomial in $z$ of degree $-a_j$.
In (\ref{e11}) the Pochhammer symbol, or ascending factorial, $(a)_n$ is given for integer $n$ by 
\[(a)_n=\frac{\g(a+n)}{\g(a)}=\left\{\begin{array}{ll} 1 & (n=0)\\a(a+1)\ldots (a+n-1) & (n\geq 1),\end{array}\right.\]
where $\g$ is the gamma function. 
Throughout we shall adopt the convention of writing the finite sequence of parameters $(a_1, \ldots , a_p)$ 
simply by $(a_p)$ and the product of $p$ Pochhammer symbols by
\[((a_p))_k\equiv (a_1)_k \ldots (a_p)_k,\]
where an empty product $p=0$ is understood to be unity.

There exist several classical summation theorems for hypergeometric series of specialized argument.
These are the theorems of Gauss, Kummer and Bailey for the ${}_2F_1$ series and Saalsch\"utz and Watson
 for the ${}_3F_2$ series; see, for example, \cite[Appendix III]{S}. Various contiguous extensions of
 these summations theorems have been obtained; see \cite{KRR, RR11} and the references therein. Recent
 work has been concerned with the extension of the above-mentioned summation theorems to higher-order
 hypergeometric series with $r$ pairs of numeratorial and denominatorial parameters differing by a set
 of positive integers $(m_r)$. One of the first results of this type is the generalized Karlsson-Minton
 summation theorem \cite{MS}, which extends the first Gauss summation theorem, given by
\bee\label{e11a}
{}_{r+2}F_{r+1}\left[\!\!\begin{array}{c}a,\, b,\\1+b,\end{array}\!\!\!\!\begin{array}{c}(f_r+m_r)\\(f_r)\end{array}\!;1\right]=\frac{\g(1-a) \g(1+b)}{\g(1+b-a)}\,\frac{(f_1-b)_{m_1} \ldots (f_r-b)_{m_r}}{(f_1)_{m_1} \ldots (f_r)_{m_r}}
\ee
provided Re$\,(-a)>m-1$, where $m:=m_1+\cdots +m_r$. When $b=-n$, where $n$ is a non-negative integer, the series on the left-hand side terminates and (\ref{e11a}) reduces to
the result originally obtained by Minton \cite{Mn} when $n\geq m$. A generalization of (\ref{e11a}) when the series terminates (an extension of the Vandermonde-Chu summation formula for the ${}_2F_1$ series) was derived by Miller \cite{M} in the form
\bee\label{e12c}
{}_{r+2}F_{r+1}\left[\!\!\begin{array}{c}-n,\,a,\\c,\end{array}\!\!\!\!\begin{array}{c}(f_r+m_r)\\(f_r)\end{array}\!;1\right]=
\frac{(c-a-m)_n}{(c)_n}\,\frac{((\xi_m+1))_n}{((\xi))_n}
\ee
for non-negative integer $n$ when it is supposed that $(c-a-m)_m\neq 0$ and $a\neq f_j$ ($1\leq j\leq r$),
where  $m$ is as defined above. The $(\xi_m)$ are the non-vanishing zeros of a certain parametric polynomial $Q_m(t)$ defined in (\ref{e31a}). The summation theorems of Gauss (second), Kummer, Bailey and Watson have been similarly extended in \cite{RP}.

The summation theorem with which we shall be concerned in this paper is Saalsch\"utz's summation theorem given by \cite[p.~49]{S}
\bee\label{e12a}
{}_3F_2\left[\!\!\begin{array}{c} -n,\, a, \, b \\c, -n-\sigma\end{array}\!; 1\right]=\frac{(c-a)_n (c-b)_n}{(c)_n (c-a-b)_n}, \qquad \sigma:=c-a-b-1
\ee
for non-negative integer $n$. An extension of this theorem to include $r$ pairs of numeratorial and denominatorial parameters differing by positive integers $(m_r)$
has recently been obtained in \cite{KRP} in the form
\bee\label{e15}
{}_{r+3}F_{r+2}\left[\!\!\begin{array}{c}-n, \,a,\, b,\\c,\,m-n-\sigma,\end{array}\!\!\!\!\!\!\begin{array}{c}(f_r+m_r)\\(f_r)\end{array}\!;1\right]
=\frac{(c-a-m)_n (c-b-m)_n}{(c)_n (c-a-b-m)_n}\,\frac{((\eta_m+1))_n}{((\eta_m))_n},
\ee
where $m$ is as above and the $(\eta_m)$ are the non-vanishing zeros of a different parametric polynomial ${\hat Q}_m(t)$ defined in (\ref{e21}).
The derivation of (\ref{e15}) relied on an Euler-type transformation for the generalized hypergeometric function
\[{}_{r+2}F_{r+1}\left[\!\!\begin{array}{c}a,\, b,\\c,\end{array}\!\!\!\!\begin{array}{c}(f_r+m_r)\\(f_r)\end{array}\!;z\right]\]
obtained in \cite{MP1, MP2}. 

It is the purpose of the present investigation to provide 
an alternative, more elementary proof of the extended Saalsch\"utz summation theorem (\ref{e15}) and also to give some illustrative examples.
Two applications are discussed involving hypergeometric series when $r$ pairs of numeratorial and denominatorial parameters differ by positive integers $(m_r)$. The first extends two transformations originally obtained by Ramanujan \cite{B} and the second
extends two reduction formulas for the Kamp\'e de F\'eriet function given by Exton \cite{E} and Cvijovi\'c and Miller \cite{CM}. 

\vspace{0.6cm}

\begin{center}
{\bf 2. \  Alternative proof of the extension of Saalsch\"utz's formula (\ref{e15})}
\end{center}
\setcounter{section}{2}
\setcounter{equation}{0}
\renewcommand{\theequation}{\arabic{section}.\arabic{equation}}
For the set of positive integers $(m_r)$ define the integer $m$ by
\[m:=m_1+\cdots +m_r.\]
Let the quantities $(\eta_m)$ be the non-vanishing zeros of the associated parametric polynomial ${\hat Q}_m(t)$ of degree $m$ given by
\bee\label{e21}
{\hat Q}_m(t)=\sum_{k=0}^m \frac{(-1)^kC_{k,r} (a)_k (b)_k (t)_k}{(c-a-m)_k (c-b-m)_k}\,G_{m,k}(t)
\ee
where
\bee\label{e22a}
G_{m,k}(t):={}_3F_2\left[\!\!\begin{array}{c}-m+k, t+k, c-a-b-m\\c-a-m+k, c-b-m+k\end{array}\!;1\right].
\ee
The coefficients $C_{k,r}$ $(0\leq k\leq m)$ are defined by
\bee\label{e22}
C_{k,r}=\frac{1}{\Lambda}\sum_{j=k}^m \sigma_{j} {\bf S}^{(k)}_j,\qquad \Lambda=(f_1)_{m_1}\ldots (f_r)_{m_r},\ \ \ C_{0,r}=1,\ C_{m,r}=1/\Lambda,
\ee
where ${\bf S}^{(k)}_j$ is the Stirling number of the second kind and the coefficients $\sigma_j$ $(0\leq j\leq m)$ are generated by
\bee\label{e22b}
(f_1+x)_{m_1} \ldots (f_r+x)_{m_r}=\sum_{j=0}^m\sigma_{j}x^j.
\ee
We note that the polynomial ${\hat Q}_m(t)$ has been normalized so that ${\hat Q}_m(0)=1$ and, for $0\leq k\leq m$, the function $G_{m,k}(t)$ is a polynomial in $t$ of degree $m-k$.
\vspace{0.2cm}

\noindent{\bf Remark 1.}\ \ \ An alternative representation of the coefficients $C_{k,r}$ is
given as the terminating hypergeometric series of unit argument \cite{MP3}
\[C_{k,r}=\frac{(-1)^k}{k!}\,{}_{r+1}F_r\left[\!\!\begin{array}{c}-k,\\{}\end{array}
\!\!\!\begin{array}{c}(f_r+m_r)\\(f_r)\end{array}\!;1\right].\]
\vspace{0.2cm}

In the statement of our theorem, we shall require the quantity $H_n$ defined by
\bee\label{e200}
H_n=\frac{((\eta_m+1))_n}{((\eta_m))_n}=\frac{(n+\eta_1)\ldots (n+\eta_m)}{\eta_1\ldots\eta_m}.
\ee
We observe that it is {\it not necessary to evaluate the zeros} $(\eta_m)$ of the associated parametric polynomial ${\hat Q}_m(t)$ to evaluate $H_n$. Since ${\hat Q}_m(t)=(-1)^m (t-\eta_1)\dots (t-\eta_m)/(\eta_1\ldots \eta_m)$, it follows immediately from (\ref{e200}) that
\[H_n={\hat Q}_m(-n),\]
where, since $(-n)_m=0$ for $m>n$,
\bee\label{e200a}
{\hat Q}_m(-n)=\sum_{k=0}^{\min \{m,n\}}\frac{(-1)^kC_{k,r} (a)_k (b)_k (-n)_k}{(c-a-m)_k (c-b-m)_k}\,G_{m,k}(-n).
\ee
Hence, to evaluate $H_n$ it suffices to determine only the associated parametric polynomial ${\hat Q}_m(t)$ and set $t=-n$.

The extension
of Saalsch\"utz's summation formula when there are $r$ pairs of numeratorial and denominatorial parameters differing by positive integers $(m_r)$ is then given by the following theorem.
\vspace{0.2cm}

\noindent{\bf Theorem 1.}\ {\it Let $(m_r)$ be a set of positive integers with $m:=m_1+\cdots +m_r$ and let $n$ denote a non-negative integer. Then, with $\sigma:=c-a-b-1$, we have \cite{KRP}
\bee\label{e20}
{}_{r+3}F_{r+2}\left[\!\!\begin{array}{c}-n, \,a,\, b,\\c,\,m-n-\sigma,\end{array}\!\!\!\!\!\!\begin{array}{c}(f_r+m_r)\\(f_r)\end{array}\!;1\right]
=\frac{(c-a-m)_n (c-b-m)_n}{(c)_n (c-a-b-m)_n}\,H_n,
\ee
where 
\[H_n:=\frac{((\eta_m+1))_n}{((\eta_m))_n}={\hat Q}_m(-n).\]
The $(\eta_m)$ are the nonvanishing zeros of the associated parametric polynomial ${\hat Q}_m(t)$ of degree $m$ defined in (\ref{e21}) and (\ref{e22a}).}

\vspace{0.4cm}

\noindent{\bf Proof:}\ \ \  To prove the summation formula (\ref{e20}), we express the ${}_{r+3}F_{r+2}(1)$ in its series form for non-negative integer $n$ as 
\[F\equiv {}_{r+3}F_{r+2}\left[\!\!\begin{array}{c}-n, \,a,\, b,\\c,\,m-n-\sigma,\end{array}\!\!\!\!\!\!\begin{array}{c}(f_r+m_r)\\(f_r)\end{array}\!;1\right]
=\sum_{s=0}^n \frac{(-n)_s (a)_s (b)_s}{(c)_s (m-n-\sigma)_s s!}\,\frac{((f_r+m_r))_s}{((f_r))_s}.\]
Making use of the fact that for non-negative integer $s$
\[\frac{(f_r+m_r)_s}{(f_r)_s}=\frac{(f_r+s)_{m_r}}{(f_r)_{m_r}},\]
we find using the definition of the coefficients $\sigma_j$ in (\ref{e22b}) and $m=m_1+\cdots +m_r$ that
\begin{eqnarray*}
\frac{((f_r+m_r))_s}{((f_r))_s}\!\!&=&\!\!\frac{1}{\Lambda} (s+f_1)_{m_1} \ldots (s+f_r)_{m_r} =\frac{1}{\Lambda}\sum_{j=0}^m\sigma_j s^j\\
&=&\!\!1+\frac{1}{\Lambda}\sum_{j=1}^m\sigma_j \sum_{k=1}^j {\bf S}_j^{(k)} s(s-1)\ldots (s-k+1)\\
&=&\!\!1+\sum_{k=1}^ms(s-1)\ldots (s-k+1) \frac{1}{\Lambda}\sum_{j=k}^m\sigma_j {\bf S}_j^{(k)}\\
&=&\!\!1+\sum_{k=1}^m C_{k,r}\,s(s-1)\ldots (s-k+1),
\end{eqnarray*}
where ${\bf S}_j^{(k)}$ is the Stirling number of the second kind and the coefficients $C_{k,r}$ and the quantity $\Lambda$ are defined in (\ref{e22}).

Then, since $(-n)_m=0$ for $m>n$ and $C_{0,r}=1$, we have
\[F=\sum_{s=0}^n \frac{(-n)_s (a)_s (b)_s}{(c)_s (m-n-\sigma)_s}
\sum_{k=0}^{\min \{m,n\}}\frac{C_{k,r}}{(s-k)!}=\sum_{k=0}^{\min \{m,n\}} C_{k,r}\sum_{s=0}^{n-k}  \frac{(-n)_{s+k} (a)_{s+k} (b)_{s+k}}{(c)_{s+k} (m-n-\sigma)_{s+k}s!}\]
upon reversing the order of summation. 
Employing the result 
\bee\label{e300}
(a)_{s+k}=(a)_k (a+k)_s, 
\ee
we obtain
\bee\label{e23}
F=\sum_{k=0}^{\min \{m,n\}} \frac{C_{k,r}\,(-n)_k (a)_k (b)_k}{(c)_k (m-n-\sigma)_k}\,
{}_3F_2\left[\!\!\begin{array}{c} -n+k,\,a+k,\,b+k\\c+k,\,m+k-n-\sigma\end{array}\!;1\right].
\ee

We now employ the contiguous Saalsch\"utz summation formula, which can be obtained from \cite[p.~539, Eq.~(85)]{PBM}, in the form 
\[{}_3F_2\left[\!\!\begin{array}{c}-n,\,a,\,b\\c,\,p-n-\sigma\end{array}\!;1\right]\hspace{10cm}\]
\[
\hspace{2cm}=\frac{(c-a-p)_n (c-b-p)_n}{(c)_n (c-a-b-p)_n}\,{}_3F_2\left[\!\!\begin{array}{c}-p,\,-n,\,c-a-b-p\\c-a-p,\,c-b-p\end{array}\!;1\right]\quad p=0, 1, 2, \ldots
\]
for non-negative integer $n$ and $\sigma:=c-a-b-1$.
Then, with $p=m-k$, $n\ra n-k$, $a\ra a+k$, $b\ra b+k$, $c\ra c+k$ and $\sigma\ra \sigma-k$, we find that the term with index $k$ in (\ref{e23}) becomes
\[\frac{C_{k,r}\,(-n)_k (a)_k (b)_k}{(c)_k (m-n-\sigma)_k}\,\frac{(c-a-m+k)_{n-k} (c-b-m+k)_{n-k}}{(c+k)_{n-k} (c-a-b-m)_{n-k}}\hspace{4cm}\]
\[\hspace{6cm}\times{}_3F_2\left[\!\!\begin{array}{c}-m+k,\,-n+k,\,c-a-b-m\\c-a-m+k,\,c-b-m+k\end{array}\!;1\right]\]
\[=\frac{(c-a-m)_n (c-b-m)_n}{(c)_n (c-a-b-m)_n}\,\frac{(-1)^kC_{k,r}\,(-n)_k (a)_k (b)_k}{(c-a-m)_k (c-b-m)_k}\,G_{m,k}(-n),\]
where we have made use of the identities (\ref{e300}) (with $s+k\ra n$) and
\[(m-n-\sigma)_k (c-a-b-m)_{n-k}=(-1)^k (c-a-b-m)_n\]
and the definition of the polynomial $G_{m,k}(t)$ in (\ref{e22a}).

Hence it follows that
\[F=\frac{(c-a-m)_n (c-b-m)_n}{(c)_n (c-a-b-m)_n} \sum_{k=0}^{\min \{m,n\}} \frac{(-1)^k C_{k,r}\,(-n)_k (a)_k (b)_k }{(c-a-m)_k (c-b-m)_k}\,G_{m,k}(-n)\]
\[=\frac{(c-a-m)_n (c-b-m)_n}{(c)_n (c-a-b-m)_n}\,{\hat Q}_m(-n)\]
by (\ref{e200a}). Since, from (\ref{e200}), 
${\hat Q}_m(-n)=((\eta_m+1))_n/((\eta_m))_n$, 
this completes the proof of the summation formula in (\ref{e20}).\ \ \ \ \ $\Box$

\vspace{0.6cm}

\begin{center}
{\bf 3. \  Examples}
\end{center}
\setcounter{section}{3}
\setcounter{equation}{0}
\renewcommand{\theequation}{\arabic{section}.\arabic{equation}}
In the case $r=1$ and $m_1=m=1$, $f_1=f$, the summation theorem (\ref{e20}) takes the form 
\bee\label{e12}
{}_4F_3\left[\!\!.\begin{array}{c}-n,\, a,\, b,\\c,\, 1-n-\sigma,\end{array}\!\!\!\begin{array}{c}f+1\\f\end{array}\!;1\right]=\frac{(c-a-1)_n (c-b-1)_n}{(c)_n (c-a-b-1)_n}\,\left(1+\frac{n}{\eta}\right)
\ee
for non-negative integer values of $n$, where
\bee\label{eeta}
\eta=\frac{(c-a-1)(c-b-1)f}{ab+(c-a-b-1)f} 
\ee
is the non-vanishing zero of the first-degree parametric polynomial 
\[{\hat Q}_1(t)=1-\frac{\{(c-a-b-1)f+ab\}t}{(c-a-1)(c-b-1)f}.\]
\vspace{0.2cm}

\noindent{\bf Remark 2.}\ \ \ When $c=1+a-b$, we have 
\bee\label{e2r}
\eta=\frac{(a-2b)f}{2f-a}.
\ee
If, in addition, $f=\fs a$ ($a\neq 2b$), we obtain from (\ref{e12}) 
\[{}_4F_3\left(\!\!\left.\begin{array}{c}-n,\, a,\, b,\ 1+\fs a\\1+a-b, 1+2b-n, \fs a\end{array}\!\right| 1\right)=\frac{(a-2b)_n (-b)_n}{(1+a-b)_n (-2b)_n}\]
which is a known result \cite[Appendix III, Eq.~(17)]{S}.
\vspace{0.2cm}

In the case $r=1$, $m_1=2$, $f_1=f$, where $C_{0,r}=1$, $C_{1,r}=2/f$ and $C_{2,r}=1/(f)_2$,
we have the quadratic parametric polynomial (with zeros $\eta_1$ and $\eta_2$) given by \cite{MP2}
\[{\hat Q}_2(t)=1-\frac{2Bt}{(c-a-2)(c-b-2)}+\frac{Ct(1+t)}{(c-a-2)_2(c-b-2)_2},\]
where
\[B:=\sigma-1+\frac{ab}{f},\qquad C:=(\sigma-1)\left(\sigma+\frac{2ab}{f}\right)+\frac{(a)_2(b)_2}{(f)_2}.\]
Hence we obtain
\[{}_4F_3\left[\!\!\begin{array}{c} -n,\,a,\,b,\\c,\,2-n-\sigma,\end{array}\!\!\!\!\!\begin{array}{c}f+2\\f\end{array}\!;1\right]=\frac{(c-a-2)_n (c-b-2)_n}{(c)_n (c-a-b-2)_n}\,\left(1+\frac{n}{\eta_1}\right)\left(1+\frac{n}{\eta_2}\right)\]
\[=\frac{(c-a-2)_n (c-b-2)_n}{(c)_n (c-a-b-2)_n}
\,\left\{1\!+\!\frac{2Bn}{(c-a-2)(c-b-2)}\!+\!\frac{Cn(n-1)}{(c-a-2)_2(c-b-2)_2}\right\}\]
for nonegative integer values of $n$.
\vspace{0.6cm}

\begin{center}
{\bf 4. \  The extension of two transformation formulas of Ramanujan}
\end{center}
\setcounter{section}{4}
\setcounter{equation}{0}
\renewcommand{\theequation}{\arabic{section}.\arabic{equation}}
For our first application of the extension of the Saalsch\"utz summation theorem in (\ref{e20}) we obtain two transformation formulas involving ${}_{r+2}F_{r+1}(x^2)$
series, when $r$ pairs of numeratorial and denominatorial parameters differ by positive integers, that generalize results originally given by Ramanujan \cite{B} in the case $r=0$.
Our results are given by the following theorem.
\vspace{0.2cm}

\noindent{\bf Theorem 2.}\ {\it Let $(m_r)$ be a set of positive integers with $m:=m_1+\cdots +m_r$.  Then, for $n$ arbitrary,
\[(1-x^2)^{-\fr} {}_{r+2}F_{r+1}\left[\!\!\begin{array}{c}-n+\fs p,\, n+\fs p,\, (f_r+m_r)\\p+\fs+m,\, (f_r)\end{array}\!;x^2\right]\hspace{6cm}\]
\bee\label{e40}
\hspace{4cm}={}_{m+2}F_{m+1}\left[\!\!\begin{array}{c}\fs+\fs p-n,\,\fs+\fs p+n,\,  (\eta_m+1)\\p+\fs+m,\, (\eta_m)\end{array}\!;x^2\right]
\ee
when $|x|<1$, where $p=0,\,1$. The $(\eta_m)$ are the non-vanishing zeros of the parametric polynomial ${\hat Q}_m(t)$ of degree $m$ in (\ref{e21}) and \ref{e22a}).}
\vspace{0.4cm}

\noindent{\bf Proof:}\ \ \  We consider the case $p=0$ and define
\begin{eqnarray*}
F\!\!&\equiv&\!\! (1-x^2)^{-\fr} {}_{r+2}F_{r+1}\left[\!\!\begin{array}{c}-n,\, n,\, (f_r+m_r)\\ \fs+m,\, (f_r)\end{array}\!;x^2\right]\\
&=&\!\!\sum_{j=0}^\infty \frac{(\fs)_j }{j!}x^{2j}\,\sum_{k=0}^\infty \frac{(-n)_k (n)_k ((f_r+m_r))_k }{(\fs+m)_k ((f_r))_k k!}x^{2k}\\
&=&\!\!\sum_{j=0}^\infty \sum_{k=0}^\infty \frac{(\fs)_j (-n)_k (n)_k ((f_r+m_r))_k }{(\fs+m)_k ((f_r))_k\, j! k!}x^{2j+2k}, \qquad (|x|<1)
\end{eqnarray*}
where we have expressed the Cauchy product as a double sum. Changing the double sum by rows to diagonal summation (see \cite[p.~58]{S}) by putting $j\rightarrow j-k$ ($0\leq k\leq j$), we find
\begin{eqnarray*}
F\!\!&=&\!\!\sum_{j=0}^\infty \sum_{k=0}^j \frac{(\fs)_{j-k} (-n)_k (n)_k ((f_r+m_r))_k}{(\fs+m)_k ((f_r))_k\,(j-k)! k!} x^{2j}\\
&=&\!\!\sum_{j=0}^\infty \frac{(\fs)_j}{j!} x^{2j} \sum_{k=0}^j\frac{(-j)_k (-n)_k (n)_k ((f_r+m_r))_k}{(\fs+m)_k (\fs+j)_k ((f_r))_k k!}\\
&=&\!\!\sum_{j=0}^\infty\frac{(\fs)_j}{j!} x^{2j}\,{}_{r+3}F_{r+2} \left[\!\!\begin{array}{c}-j,\,-n,\,n,\,(f_r+m_r)\\\fs-j,\ \fs+m,\,(f_r)\end{array}\!;1\right],
\end{eqnarray*}
where we have made use of the identities
\bee\label{eid}
(\fs)_{j-k}=\frac{(-1)^k (\fs)_k}{(\fs-j)_k},\qquad \frac{1}{(j-k)!}=\frac{(-1)^k (-j)_k}{j!}.\ee

If we now identify the parameters $a$, $b$ and $c$ in (\ref{e20}) with $n$, $-n$ and $\fs+m$ respectively, then we can apply the extension of Saalsch\"utz's summation formula in Theorem 1 to obtain
\[F=\sum_{j=0}^\infty \frac{(\fs-n)_j (\fs+n)_j ((\eta_m+1))_j}{(\fs+m)_j ((\eta_m))_j\,j!} x^{2j}={}_{m+2}F_{m+1}\left[\!\!\begin{array}{c}\vspace{0.1cm}

\fs-n,\,\fs+n,\,(\eta_m+1)\\\fs+m,\,(\eta_m)\end{array}\!;x^2\right],\]
thereby establishing the result when $p=0$. The proof of the case with $p=1$ is similar and consequently will be omitted. \ \ \ \ \ $\Box$
\vspace{0.4cm}

In the case $r=0$ ($m=0$), (\ref{e40}) reduces to the two identities for $n$ arbitrary
\bee\label{e41a}
(1-x^2)^{-\fr}\,{}_2F_1\left[\!\!\begin{array}{c} -n,\,n\\ \fs\end{array}\!;x^2\right]={}_2F_1\left[\!\!\begin{array}{c}\vspace{0.1cm}

\fs-n,\,\fs+n\\\fs\end{array}\!;x^2\right]
\ee
and
\bee\label{e41b}
(1-x^2)^{-\fr}\,{}_2F_1\left[\!\!\begin{array}{c} \fs-n,\,\fs+n\\ \f{3}{2}\end{array}\!;x^2\right]={}_2F_1\left[\!\!\begin{array}{c}\vspace{0.1cm}

1-n,\,1+n\\\f{3}{2}\end{array}\!;x^2\right]
\ee
obtained by Ramanujan \cite[p.~99, 35(iii)]{B} and \cite{CR}, respectively.

When $r=1$, $m_1=m=1$, $f_1=f$, we obtain from (\ref{e40}) with $p=0$
\bee\label{e41}
(1-x^2)^{-\fr}\,{}_3F_2\left[\!\!\begin{array}{c}\vspace{0.1cm}

 -n,\,n,\,f+1\\ \f{3}{2},\,f\end{array}\!;x^2\right]={}_3F_2\left[\!\!\begin{array}{c}\vspace{0.1cm}

\fs-n,\,\fs+n,\,\eta+1\\ \f{3}{2},\,\eta\end{array}\!;x^2\right],
\ee
where, from (\ref{eeta}), 
\[\eta=\frac{(n^2-\f{1}{4})f}{n^2-\fs f};\]
and when $p=1$
\bee\label{e42}
(1-x^2)^{-\fr}\,{}_3F_2\left[\!\!\begin{array}{c}\vspace{0.1cm}

\fs-n,\,\fs+n,\,f+1\\ \f{5}{2},\,f\end{array}\!;x^2\right]={}_3F_2\left[\!\!\begin{array}{c}\vspace{0.1cm}

1-n,\,1+n,\,\eta+1\\ \f{5}{2},\,\eta\end{array}\!;x^2\right],
\ee
where, from (\ref{eeta}), 
\[\eta=\frac{(n^2-1)f}{n^2-\fs f-\f{1}{4}}\] 
both for $n$ arbitrary.
When $f=\fs$ and $f=\f{3}{2}$, we remark that (\ref{e41}) and (\ref{e42}) correctly reduce to (\ref{e41a}) and (\ref{e41b}), respectively.
The results in (\ref{e41}) and (\ref{e42}) have been obtained recently by different means in \cite{CRC}.

\vspace{0.6cm}

\begin{center}
{\bf 5. \  Two reduction formulas for the Kamp\'e de F\'eriet function}
\end{center}
\setcounter{section}{5}
\setcounter{equation}{0}
\renewcommand{\theequation}{\arabic{section}.\arabic{equation}}
The Kamp\'e de F\'eriet function is a hypergeometric  function of two variables defined by
\bee\label{e30}
F^{p:\,r;\,s}_{\,q:\, t;\,u}\left[\left.\!\!\begin{array}{c}(\alpha_p):\,(a_r);\,(b_s)\\(\beta_q):\,(c_t);\,(d_u)\end{array}\!\right|x,\ y\right]=
\sum_{m=0}^\infty \sum_{n=0}^\infty \frac{((\alpha_p))_{m+n}}{((\beta_q))_{m+n}}\,\frac{((a_r))_m ((b_s))_n}{((c_t))_m ((d_u))_n}\,\frac{x^my^n}{m!\,n!},
\ee
where $p$, $q$, $r$, $s$, $t$, $u$ are nonnegative integers that correspond to the number of elements in the parameter sets $(\alpha_p)$, $(\beta_q)$, $(a_r)$, $(b_s)$, $(c_t)$ and $(d_u)$, respectively; for an introduction to this function, see \cite[pp.~63--64]{SM}. We also have the easily established result \cite[Eq.~(6)]{E}
\[F^{p:\,r;\,s}_{\,q:\, t;\,u}\left[\left.\!\!\begin{array}{c}(\alpha_p):\,(a_r);\,(b_s)\\(\beta_q):\,(c_t);\,(d_u)\end{array}\!\right|x,\ y\right]\hspace{8cm}\]
\bee\label{e31}
\hspace{2cm}=\sum_{n=0}^\infty\frac{((\alpha_p))_n ((b_s))_n}{((\beta_q))_n ((d_u))_n}\,\frac{y^n}{n!}\,{}_{r+u+1}F_{s+t}\left[\!\!\begin{array}{c}-n,\,(a_r),\,(1-d_u-n)\\(c_t),\,(1-b_s-n)\end{array}\!;(-1)^{s-u+1} \frac{x}{y}\right].
\ee

Reduction formulas represent the Kamp\'e de F\'eriet function as a generalized hypergeometric function of lower order and of a single variable. The identification of such reductions is of considerable utility in the application of these functions; a compilation can be found in \cite[pp.~28--32]{SK}. In this section we shall be concerned with reduction formulas 
for the Kamp\'e de F\'eriet function when one set of numeratorial and denominatorial parameters differs by positive integers $(m_r)$. One of the first results of this type was obtained by Miller \cite{M} in the form\footnote{In \cite{M}, Miller gave the case $(m_r)=1$ but his arguments are easily extended to the case of positive integers $(m_r)$.}
\[F^{p:\,r+1;\,0}_{\,q:\, r+1;\,0}\left[\left.\!\!\begin{array}{c}(\alpha_p):\\(\beta_q):\end{array}\!\!\!\begin{array}{c}a,\,(f_r+m_r);\\b,\, (f_r);\end{array}\!\!\!\begin{array}{c}\rule{0.3cm}{0.02cm}\\ \rule{0.3cm}{0.02cm}\end{array}\!\right|-x,\ x\right]
={}_{p+m+1}F_{q+m+1}\left[\!\!\begin{array}{c}(\alpha_p),\\(\beta_q),\end{array}\!\!\!\!\begin{array}{c}b-a-m,\,(\xi_m+1)\\b,\,(\xi_m)\end{array}\!;x\right],\]
where the horizontal line indicates an empty parameter sequence. The $(\xi_m)$ are the nonvanishing zeros of the associated parametric polynomial of degree $m=m_1+\cdots +m_r$ given by
\bee\label{e31a}
Q_m(t)=\frac{1}{(\lambda)_m}\sum_{k=0}^m (b)_k C_{k,r} (t)_k (\lambda-t)_{m-k}
\ee
which is normalized so that $Q_m(0)=1$, where $\lambda:=b-a-m$ and the coefficients $C_{k,r}$ are defined in (\ref{e22}). In the case $r=1$, $m_1=m=1$, $f_1=f$, we have
\[Q_1(t)=1+\frac{(b-f)t}{(c-b-1)f},\]
with the nonvanishing zero $\xi_1=\xi$ (provided $c-b-1\neq 0$) given by
\bee\label{e31b}
\xi=\frac{(c-b-1)f}{f-b}.
\ee 

Here we shall exploit the result in (\ref{e20}) and the extension of the Vandermonde-Chu summation formula in (\ref{e12c}) to obtain two new reduction formulas. 
\vspace{0.4cm}

\noindent{\bf 5.1 \ First reduction formula}
\vspace{0.2cm}

\noindent From (\ref{e31}), with $x=y$ we obtain
\[F^{p:\,r+2;\,1}_{\,q:\, r+1;\,0}\left[\left.\!\!\begin{array}{c}(\alpha_p):\\(\beta_q):\end{array}\!\!\!\begin{array}{c}a,\,b,\ (f_r+m_r);\\c,\,(f_r);\end{array}\!\!\!\begin{array}{c}c-a-b-m\\ \rule{0.3cm}{0.02cm}\end{array}\!\right|x,\ x\right]\hspace{8cm}\]
\begin{eqnarray}
&=&\!\!\sum_{n=0}^\infty \frac{((\alpha_p))_n}{((\beta_q))_n}\,\frac{(c-a-b-m)_n}{n!} x^n\,{}_{r+3}F_{r+2}\left[\!\!\begin{array}{c}-n,\,a,\,b,\\c,\, m-n-\sigma,\end{array}
\!\!\!\!\!\!\!\begin{array}{c}(f_r+m_r)\\(f_r)\end{array}\!;1\right]\nonumber\\
&=&\!\!\sum_{n=0}^\infty \frac{((\alpha_p))_n}{((\beta_q))_n}\,\frac{(c-a-m)_n (c-b-m)_n}{(c)_n}\,\frac{((\eta_m+1))_n}{((\eta_m))_n}\,\frac{x^n}{n!}\nonumber\\
&=&\!\!{}_{p+m+2}F_{q+m+1}\left[\!\!\begin{array}{c}(\alpha_p),\\(\beta_q),\end{array}\!\!\!\begin{array}{c}c-a-m, c-b-m,\\c,\end{array}\!\!\!\begin{array}{c}(\eta_m+1)\\(\eta_m)\end{array}\!;x\right]\label{e33}
\end{eqnarray}
upon application of (\ref{e20}), where we recall that $\sigma:=c-a-b-1$ and the $(\eta_m)$ are the nonvanishing zeros of the parametric polynomial ${\hat Q}_m(t)$ defined in (\ref{e21}).

In the case $r=1$, $m_1=m=1$, $f_1=f$ we obtain the reduction formula
\[F^{p:\,3;\,1}_{\,q:\, 2;\,0}\left[\left.\!\!\begin{array}{c}(\alpha_p):\\(\beta_q):\end{array}\!\!\!\begin{array}{c}\,a,\,b,\,f+1;\\c,\, f;\end{array}\!\!\!\begin{array}{c}c-a-b-1\\ \rule{0.3cm}{0.02cm}\end{array}\!\right|x,\ x\right]\hspace{6cm}\]
\bee\label{e302}
\hspace{5cm}={}_{p+3}F_{q+2}\left[\!\!\begin{array}{c}(\alpha_p),\\(\beta_q),\end{array}\!\!\!\begin{array}{c} c-a-1, c-b-1, \,\eta+1\\c,\,\eta\end{array}\!;x\right],\label{e32}
\ee
where $\eta$ is defined in (\ref{eeta}).
\vspace{0.2cm}

\noindent{\bf Remark 3:}\ \ \ If we take $c=1+a-b$ in (\ref{e302}), then $\eta$ is given by  (\ref{e2r}).
The special case  $f=\fs a$  then corresponds to $\eta=\infty$ ($a\neq 2b$) and we obtain Exton's reduction formula \cite[Eq.~(11)]{E}
\[F^{p:\,3;\,1}_{\,q:\, 2;\,0}\left[\left.\!\!\begin{array}{c}(\alpha_p):\\(\beta_q):\end{array}\!\!\!\begin{array}{c}\,a,\,b,\,1+\fs a;\\1+a-b,\ \fs a;\end{array}\!\!\!\begin{array}{c}-2b\\ \rule{0.3cm}{0.02cm}\end{array}\!\right|x,\ x\right]
={}_{p+2}F_{q+1}\left[\!\!\begin{array}{c}(\alpha_p),\\(\beta_q),\end{array}\!\!\!\begin{array}{c} a-2b,\, -b\\1+a-b\end{array}\!;x\right].\]
\vspace{0.4cm}

\noindent{\bf 5.2 \ Second reduction formula}
\vspace{0.2cm}

\noindent We consider another reduction formula when $p=q=1$ for the function
\[F(x,x)\equiv F^{1:\,r+2;\,1}_{\,1:\, r+1;\,0}\left[\left.\!\!\begin{array}{c}\alpha:\\ \beta:\end{array}\!\!\!\begin{array}{c}a,\,\beta-d,\, (f_r+m_r);\\c,\,(f_r);\end{array}\!\!\!\begin{array}{c}d\\ \rule{0.3cm}{0.02cm}\end{array}\!\right|x,\ x\right]\]
thereby generalizing a result obtained by Cvijovi\'c and Miller in \cite{CM} in the case $r=0$. Our derivation follows closely that given by these authors. 

Making use of the identity (\ref{e300}) we find when $|x|<1$
\begin{eqnarray*}
F(x, x)\!\!&=&\!\!\sum_{k=0}^\infty\sum_{n=0}^\infty\frac{(\alpha)_{n+k}}{(\beta)_{n+k}}\,\frac{(a)_k (\beta-d)_k (d)_n}{(c)_k}\,\frac{((f_r+m_r))_k}{((f_r))_k}\,\frac{x^{n+k}}{n!\,k!}\\
&=&\!\!\sum_{k=0}^\infty \frac{(\alpha)_k (a)_k (\beta-d)_k}{(\beta)_k (c)_k}\,\frac{((f_r+m_r))_k}{((f_r))_k} \frac{x^k}{k!}\,{}_2F_1\left[\!\!\begin{array}{c} \alpha+k,\,d\\\beta+k\end{array}\!;x\right].
\end{eqnarray*}
Now applying Euler's first transformation \cite[p.~31]{S}
\[{}_2F_1\left[\!\!\begin{array}{c}a,\,b\\c\end{array}\!;x\right]=(1-x)^{-a}\,
{}_2F_1\left[\!\!\begin{array}{c}a,\,c-b\\c\end{array}\!;\frac{x}{x-1}\right],\]
we obtain
\begin{eqnarray*}
F(x,&&\hspace{-0.9cm}x)\\
&=&\!\!(1-x)^{-\alpha}\sum_{k=0}^\infty 
\frac{(\alpha)_k (a)_k (\beta-d)_k}{(\beta)_k (c)_k k!}\,\frac{((f_r+m_r))_k}{((f_r))_k} \left(\!\frac{x}{1-x}\!\right)^k\,{}_2F_1\left[\!\!\begin{array}{c} \alpha+k,\,\beta-d+k\\\beta+k\end{array}\!;\frac{x}{x-1}\right]\\
&=&\!\!(1-x)^{-\alpha}\sum_{k=0}^\infty\sum_{n=0}^\infty
\frac{(-1)^k(\alpha)_k (a)_k (\beta-d)_k}{(\beta)_k (c)_k k!}\,\frac{((f_r+m_r))_k}{((f_r))_k} \left(\!\frac{x}{x-1}\!\right)^{n+k}\!\frac{(\alpha+k)_n (\beta-d+k)_n}{(\beta+k)_n\,n!}\\
&=&\!\!(1-x)^{-\alpha}\sum_{k=0}^\infty\sum_{n=0}^\infty
\frac{(\alpha)_{n+k}(\beta-d)_{n+k}}{(\beta)_{n+k}}\,\frac{(-1)^k (a)_k((f_r+m_r))_k}{(c)_k((f_r))_k\,n!\,k!} \left(\!\frac{x}{x-1}\!\right)^{n+k}
\end{eqnarray*}
valid when $|x|<1$ and $|x/(x-1)|<1$; that is in the domain $|x|<1$, $\Re (x)<\fs$.

Making the change of summation index $n\ra n-k$, reversing the order of summation
and using the second identity in (\ref{eid}), we then find
\begin{eqnarray}
F(x, x)\!\!&=&\!\!(1-x)^{-\alpha}\sum_{n=0}^\infty
\frac{(\alpha)_{n}(\beta-d)_{n}}{(\beta)_{n} n!}\left(\!\frac{x}{x-1}\!\right)^{n}
\sum_{k=0}^n\frac{(-n)_k (a)_k ((f_r+m_r))_k}{(c)_k ((f_r))_k\,k!}\nonumber\\
&=&\!\!(1-x)^{-\alpha}\sum_{n=0}^\infty
\frac{(\alpha)_{n}(\beta-d)_{n}}{(\beta)_{n} n!}\left(\!\frac{x}{x-1}\!\right)^{n}
{}_{r+2}F_{r+1}\left[\!\!\begin{array}{c}-n,\,a,\\c,\end{array}\!\!\!\begin{array}{c}(f_r+m_r)\\(f_r)\end{array}\!;1\right]\nonumber\\
&=&\!\!(1-x)^{-\alpha}\,{}_{m+3}F_{m+2}\left[\!\!\begin{array}{c}\alpha,\,\beta-d,\,c-a-m,\\ \beta,\, c,\end{array}\!\!\!\begin{array}{c}(\xi_m+1)\\(\xi_m)\end{array}\!;\frac{x}{x-1}\right]
\end{eqnarray}
upon use of the extension of the Vandermonde-Chu summation theorem in (\ref{e12c}), where the $(\xi_m)$ are the nonvanishing zeros of the associated parametric polynomial $Q_m(t)$ defined in (\ref{e31a}).

When $r=0$ ($m=0$) we recover the Cvijovi\'c-Miller result \cite{CM}  given by 
\[F^{1:\,2;\,1}_{\,1:\, 1;\,0}\left[\left.\!\!\begin{array}{c}\alpha:\\ \beta:\end{array}\!\!\!\begin{array}{c}a,\,\beta-d;\\c;\end{array}\!\!\!\begin{array}{c}d\\ \rule{0.3cm}{0.02cm}\end{array}\!\right|x,\ x\right]
=(1-x)^{-\alpha}{}_3F_2\left[\!\!\begin{array}{c}\alpha,\,\beta-d,\,c-a\\ \beta,\,c
\end{array}\!;\frac{x}{x-1}\right].\]
When $r=1$, $m_1=m=1$, $f_1=f$, we have the reduction formula
\[F^{1:\,3;\,1}_{\,1:\, 2;\,0}\left[\left.\!\!\begin{array}{c}\alpha:\\ \beta:\end{array}\!\!\!\begin{array}{c}a,\,\beta-d,\, f+1;\\c,\,f;\end{array}\!\!\!\begin{array}{c}d\\ \rule{0.3cm}{0.02cm}\end{array}\!\right|x,\ x\right]\hspace{6cm}\]
\bee
\hspace{5cm}=
(1-x)^{-\alpha}{}_4F_3\left[\!\!\begin{array}{c}\alpha,\,\beta-d,\,c-a-1,\\ \beta,\,c,\end{array}\!\!\!\!\begin{array}{c} \xi+1\\ \xi\end{array}\!;\frac{x}{x-1}\right],
\ee
where, from (\ref{e31b}) with $b$ replaced by $a$, 
\[\xi=\frac{(c-a-1)f}{f-a}.\]

\vspace{0.6cm}

\noindent{\bf Acknowledgement:}\ \ \ Y. S. Kim acknowledges the support of the Wonkwang University Research Fund (2013).

\end{document}